\documentclass[twocolumn]{article}
\usepackage{graphicx}
\usepackage{epsfig}
\usepackage{url}
\usepackage{amssymb}
\usepackage{amsmath}
\usepackage{epstopdf}
\usepackage{caption}
\usepackage{subcaption}
\usepackage{graphicx}
\usepackage{dcolumn}
\usepackage{bm}
\usepackage{hyperref}
\usepackage{color}

\usepackage{amsthm}
\theoremstyle{plain}
\newtheorem{thm}{Theorem}

\title{Evaluating the Usefulness of Paratransgenesis for Malaria Control}
\author{Bhushan Kotnis\footnote{Indian Institute of Science, bkotnis@dese.iisc.ernet.in} and Joy Kuri\footnote{Indian Institute of Science}}

\begin{document}
\maketitle

\begin{abstract}
Malaria is a serious global health problem which is especially devastating to the developing world. Most malaria control programs use insecticides for controlling mosquito populations. Large scale usage of these insecticides exerts massive selection pressure on mosquitoes resulting in insecticide resistant mosquito breeds. Thus, developing alternative strategies are crucial for sustainable malaria control. Here, we explore the usefulness of an alternative strategy, {paratransgenesis:} the introduction of genetically engineered plasmodium killing bacteria inside the mosquito gut. The genetically modified bacterial culture is housed in cotton balls dipped in a sugar solution (sugar bait) and they enter a mosquito's midgut when it drinks from a sugar bait. We study scenarios where vectors and hosts mix homogeneously as well as heterogeneously and calculate the amount of baits required to prevent a malaria outbreak. Given the baits are attractive, we show that the basic reproductive number drops rapidly with the increase in bait density. Furthermore, we propose a targeted bait distribution strategy for minimizing the reproductive number for the heterogeneous case. Our results can prove to be useful for designing future experiments and field trials of alternative malaria control mechanisms and they also have implications on the development of malaria control programs.  
\end{abstract}

\section{Introduction \label{sec:Introduction}}
The spread of malaria is a serious health concern worldwide. Malaria alone is responsible for { about six hundred thousand to one million deaths per year \cite{Malaria2015}}, and for infecting 300-500 million people every year. Malaria is particularly devastating for developing countries, resulting in a 1.3 \% loss of annual G.D.P growth \cite{Sachs2002}. The primary causative agent of malaria is the parasite \emph{Plasmodium} which depends on the vector, the female \emph{Anopheles} mosquito, for completing its life cycle. { The parasite houses itself inside the salivary glands of the mosquito for gaining entry inside the human host's blood stream when the mosquito takes a blood meal.} Conversely, hosts infected by \emph{Plasmodium} can transfer the parasite to a mosquito when the latter takes a blood meal. Malaria cannot spread without mosquitoes; hence controlling the vector population, mosquito bites, or interfering in the ability of mosquitoes to house \emph{Plasmodium} can limit the spread of malaria. Currently, no vaccine exists to prevent malaria, and hence efforts to control malaria are primarily based on vector control \cite{Enayati2010}. 

Insecticides have been very useful for controlling mosquitoes \cite{Enayati2010,Raghavendra2011}, and therefore control measures heavily depend on selective indoor residual spraying (IRS) and insecticide treated nets (ITN). The industrial scale usage of insecticides for controlling mosquito populations exerts massive selection pressure on mosquitoes. 
This evolutionary stress has resulted in insecticide resistant mosquito breeds \cite{Raghavendra2011,Read2009}. Although development of new insecticides can address these problems, it is just a matter of time before mosquitoes develop resistance to these new insecticides. Another important problem with vector control is its continual administration \cite{Wang2013}. An interruption in insecticide treatment will cause the mosquito population to rebound to pre-treatment levels. There is a great need to develop novel approaches for long term sustainable control of malaria. The malERA consulting group, in a recent report, has stressed that developing innovative strategies is crucial for sustainable vector control on a global scale \cite{Baum2011}. To fill this need, we explore one such innovative mechanism.

 The introduction of genetically modified viruses or bacteria, which can thrive in the mosquito's midgut and kill the parasite is termed as paratransgenesis \cite{Wang2012}.{ This is different than \emph{transgenesis}, i.e., modifying the mosquito genetically to impair their malaria parasite carrying capacity and releasing them in the wild to replace the wild type mosquitoes. Mosquito genes can be modified in the lab so that the mosquito produces proteins that either inhibit parasite reproduction, or kill the parasite. A variety of lab experiments \cite{Ito2002,Moreira2002} have shown that such genetic modification can reduce \emph{Plasmodium} transmission. However, simply modifying the mosquitoes and releasing them in the wild may not be enough to prevent the spread of malaria \cite{Wang2013}. To this end a variety of gene drive mechanisms using viruses (\emph{Anopheles.gambiae densonucleosis virus}  \cite{Ren2008}), and bacteria (\emph{Wolbachia} \cite{Bian2013} ) have been proposed to modify genes of wild type mosquitoes. Here, we focus on an alternative strategy - instead of modifying mosquito genes, genetically modifying bacteria are introduced inside the mosquito midgut to kill the parasite.   }
 \par
  In a recent study \cite{Wang2012}, researchers introduced a genetically modified common mosquito symbiotic bacterium \emph{Pantoea agglomerans} inside the mosquito's midgut. The genetically modified \emph{P. agglomerans} produced a variety of anti-plasmodium molecules which resulted  in {up to 98\% } reduction in the \emph{Plasmodium falciparum} population inside the midgut. Cotton balls dipped in a sugar solution containing the bacterial culture acted as baits. Thus, bacteria were introduced inside the mosquitoes when they took a sugar meal from the bait.
  \par
 The implementation of such novel techniques is crucial for sustainable control and eradication of malaria. A quantitative study of the effectiveness of these novel strategies is needed before a full scale implementation. Here, we use mathematical models to do the same. The classical Ross-McDonald model proposed in the late 1950s has exerted a large influence in modeling malaria as most models proposed from 1970-2010 are not very different from the Ross-McDonald model \cite{Smith2004a,Reiner2013,Smith2012}. Mathematical and data driven models have been used extensively for estimating the basic reproductive number \cite{Mandal2011}, which is the key quantity that determines the chance of an outbreak. Through the insights gained by analyzing various mathematical and computational models (see \cite{Mandal2011} and \cite{Reiner2013} for a recent review), researchers have formulated and evaluated various control strategies which use : insecticides \cite{Read2009b,Yakob2011,Hancock2009}, insecticide treated nets \cite{Blayneh2009,Kiware2012}, larval control \cite{Luz2011}, odor baited traps \cite{Okumu2010}, entomopathogenic fungi \cite{Lynch2012}, \emph{wolbachia} \cite{Hancock2011} and genetic modification of mosquitoes \cite{Boete2002,Atkinson2007}. However, to the best of our knowledge, no mathematical model has been formulated to evaluate the usefulness of paratransgenesis in malaria control. 
 \par
 Here we aim to study the usefulness of introducing genetically modified \emph{P. agglomerans} inside mosquitoes through sugar baits for controlling malaria, with a focus on quantifying the amount of baits required to prevent an outbreak. { Although, like insecticide spraying, sugar baits need to be replenished at regular intervals, replenishing sugar-baits is more economical and environmentally friendly than the continual usage of insecticides. Furthermore, since this strategy does not involve the killing of mosquitoes, its usage would not put evolutionary pressure on mosquitoes. An objection may be raised that the usage of anti-plasmodium  molecules may result in plasmodium which are resistant to these molecules. This  problem can be avoided by using a variety of anti plasmodium molecules \cite{Wang2013}.}
 \par
 In the past, most studies have assumed homogeneous mixing between mosquitoes and hosts \cite{Reiner2013}. However, in reality vectors and hosts may not be well mixed \cite{Woolhouse1997,Smith2005}. In this paper, we not only study the homogeneous mixing case, but also the \emph{heterogeneous} case. { Studies \cite{Smith2007,Bousema2012,Smith2010} suggest that heterogeneous mixing between the vectors and hosts may increase the basic reproductive number.} Therefore, for the heterogeneous mixing scenario we propose an optimal targeted bait allocation strategy for reducing the reproductive number. We use the \emph{Susceptible Infected Susceptible} (SIS) model for human hosts and \emph{Susceptible Infectious Removed Susceptible} (SEIRS)  model with delay for vectors. This is one of the first models which incorporates a removed compartment for mosquitoes motivated by paratransgenesis.

 Our contributions are summarized as follows:
  \begin{itemize}
 \item{{ The reproductive number} is calculated and a stability analysis is performed for the homogeneous mixing case.}
 \item{The conditions required for a malaria outbreak  is calculated for the heterogeneous mixing case.}
 \item{An optimal targeted bait allocation strategy is proposed for the heterogeneous mixing case.}
 \item{We analytically show that the reproductive number is inversely proportional to the square of effective baits for high bait attractiveness. We also discover that improving the attractiveness of sugar baits is more fruitful than increasing the efficacy of paratransgenesis.}
 \end{itemize}        
 
 { The article is organized as follows:} the model is introduced in Section \ref{sec:Model},  rigorously analyzed in Section \ref{sec:Analysis}, results detailed in Section \ref{sec:Numerical}, and the implications and interpretations of the results are discussed in Section \ref{sec:Discussion}.

 \section{Model \label{sec:Model}}
We use a compartmental SIS model for human hosts and SIERS  model with delay for mosquitoes. Human hosts can either be in the susceptible state, or the infected state. A susceptible human host can be infected by the disease, when a mosquito carrying the parasite bites the host. We assume that due to the availability of effective malarial medications, infected individuals recover at a constant rate and can be re-infected, but they do not die from the disease. After contracting the infection, mosquitoes become infectious after a fixed time duration (incubation time) which depends on the parasite. 
\par 
 When a wild type mosquito takes a blood meal from a host infected with the \emph{plasmodium} parasite, the parasite enters the mosquito's midgut in the form of \emph{gametocytes}, which reproduce, eventually producing \emph{sporozoites}, which then invade the mosquito's salivary glands. If the mosquito's midgut contains the genetically modified bacteria then the parasite will be killed by the bacteria. However, if the sporozoites have already invaded the salivary gland then the bacteria are unable to kill them. This is because the bacteria are housed in the midgut and the anti-\emph{plasmodium} molecules secreted by them cannot enter the salivary glands.  There is a time lag (in days) between the introduction of gametocytes inside the mosquito's midgut and sporozoite invasion of its salivary glands.  Based on these observations \cite{Baton2005,Wang2012}, we divide the mosquito population into four classes: susceptible, { exposed}, infectious, recovered. We incorporate the time lag (exposed class) using a fixed time delay $\tau$. Mosquitoes who do not have both: malaria parasites in their body and the genetically modified bacteria in their midgut belong to the \emph{susceptible} class. {Mosquitoes with gametocytes in midgut but no sporozoites in salivary glands and no bacteria in the midgut belong to the \emph{exposed} class. }Mosquitoes whose salivary glands are invaded by the parasites belong to the \emph{infectious} class. 
 \par
    \begin{figure}
    \centering
           \includegraphics[width=0.44\textwidth]{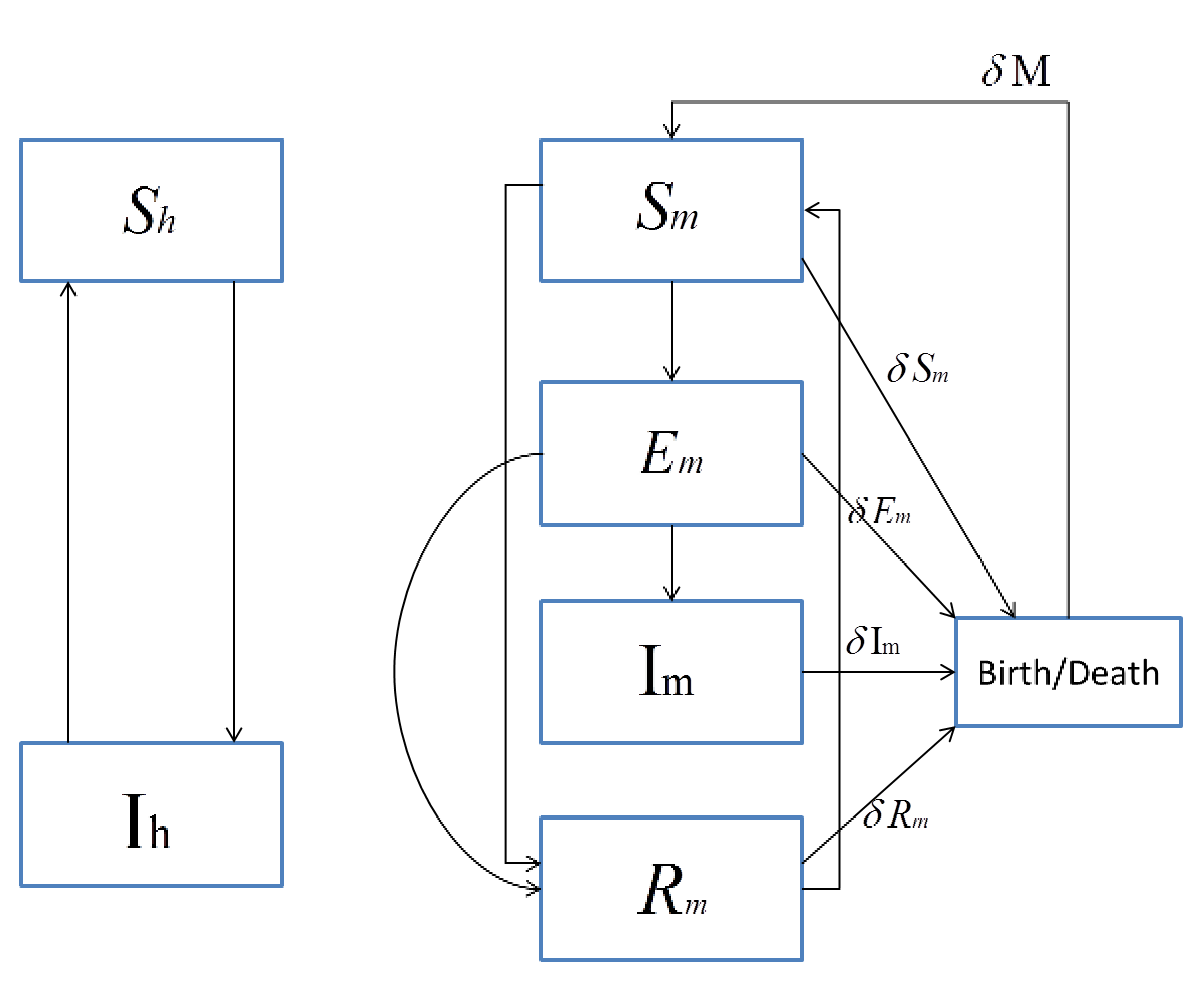}
           \caption{Block diagram of the model. Arrows represent transitions.}
           \label{fig:model_diagram}
           \end{figure}
  The effect of paratransgenesis is studied by explicitly including the \emph{removed} compartment in the model. Genetically modified \emph{P. Aggloromerans} culture is mixed with a sugar solution to act as sugar baits. Such  baits are positioned in places well within the range of mosquitoes. If a susceptible mosquito ingests the genetically engineered bacteria then it can no longer carry the parasite in its gut and enters the removed state. {Although authors in \cite{Wang2012} do not carry out experiments wherein the bacteria is introduced \emph{after} the mosquitoes take a blood meal, we believe it is reasonable to assume a similar effect if such an experiment was carried out. This is because, as the authors report, that the population of \emph{Pantoea agglomerans} increases rapidly and peaks just after 2 days. A mosquito that has ingested a blood meal, and ingests the bacteria during the incubation period, would stop carrying the parasite, because the bacteria would proliferate rapidly and release anti-plasmodium molecules that kill the oocysts in the midgut (which is also the place where the bacteria reside).   Hence we assume than a mosquito in exposed state that ingests the bacteria enters the removed state. }
  
   Drinking from a sugar bait will have no effect on infectious mosquitoes as the sporozoites have already invaded their salivary glands. The effect of  genetically modified \emph{P. Aggloromerans} bacteria may not be permanent \cite{Wang2012} and the mosquito may revert to the susceptible state. The life cycle of mosquitoes is much shorter than that of human hosts and hence must be accounted for in the model. We assume a fixed birth and death rate for the vectors, independent of their state. The birth rate is assumed to be the same as the death rate, leading to constant number of vectors. Vertical transfer of the parasite or the genetically engineered bacteria is not possible, and hence newly born mosquitoes belong to the susceptible state. Fig. \ref{fig:model_diagram} provides a graphical representation of the model.

  \begin{table}
   \small\addtolength{\tabcolsep}{-5pt}
   \caption{Definitions of symbols.}
   \label{table:Symbol Table}
 
  \begin{tabular}{l  l}
\hline\noalign{\smallskip}
  Symbol   & Description \\
 \noalign{\smallskip}\hline\noalign{\smallskip}
  $N$      & Number of human hosts. \\
  $M$      & Number of mosquitoes. \\
  $B$      & Number of sugar baits.\\
  $S_h$    & Mean number of susceptible hosts at time $t$. \\
  $I_h$    & Mean number of infected hosts at time $t$. \\
  $S_m$    & Mean number of susceptible mosquitoes at time $t$. \\
  $I_m$    & Mean number of infected mosquitoes at time $t$. \\
  $R_m$    & Mean number of removed mosquitoes at time $t$. \\
  $a$      & Biting rate of each mosquito. \\
  $b$      & Probability of infection for human host. \\ 
  $c$      & Probability of infection for mosquitoes. \\
  $m$      & Mosquito density ($M/N$).\\
  $x$      & Bait density ($B/N$).\\
  $p$     & Mosquito's preference for blood meal. \\
  $q$     & Mosquito's preference for sugar meal ($1-p$). \\
  $\tau$ & Incubation period inside mosquitoes. \\
  $\gamma$ & Efficacy of paratransgenesis.\\
  $1/\theta$ & Persistence of paratransgenesis. \\
  $\mu$    & Rate of recovery for an infected human host.\\
  $\delta$ & Birth/death rate per mosquito. \\
\noalign{\smallskip}\hline
  \end{tabular}

        \end{table}
    
Let $N$ and $M$ be the number of human hosts and mosquitoes respectively. Human populations can be homogeneous or heterogeneous with respect to their ability to attract mosquitoes; we address both the cases separately. Mosquitoes are assumed to be homogeneous in terms of their biting rate and their ability to house the genetically modified bacterium, as \emph{P. Aggloromerans} is symbiotic to most species \cite{Wang2012}. Let the mean number of healthy human hosts be $S_h$, mean infected hosts $I_h$, average number of susceptible mosquitoes be $S_m$, {mean number of exposed mosquitoes be $E_m$}, mean infectious mosquitoes be $I_m$, and mean number of removed mosquitoes be $R_m$. Let $a$ be the biting rate of each mosquito, i.e, number of bites per unit time, let $b$ be the chance that a human catches an infection due to a bite from an infected mosquito, and let $c$ be the probability that a mosquito ingests the parasite from infected human hosts. Let $B$ be the number of sugar baits, and let $x$ be the bait density, i.e., number of baits per human host: $x = B/N$. Similarly, let $m$ be the mosquito density: $m = M/N$. Let $\gamma$ be the probability that the introduced bacterium will kill all the \emph{P. falicparum} sporozoites, and let $1/\theta$ denote the persistence of paratransgenesis, i.e., the mean duration after which the bacteria become ineffective in killing the parasite. Let $\delta$ be the birth rate and the death rate (births per unit time) for mosquitoes, and let $\mu$ be the recovery rate (recoveries per unit time) for infected humans.  
  Table \ref{table:Symbol Table} provides a description of the symbols. 
  
 
\section{Analysis \label{sec:Analysis}}
Since the underlying stochastic process is Poisson, mean field approximation can be used to analyze the model. Such an approximation allows us translate the model to a system of  differential equations. The analysis for the homogeneous and the heterogeneous mixing scenarios between vectors and hosts is carried out in Section \ref{subsec:Homogeneous} and \ref{subsec:Heterogeneous} respectively. In Section \ref{subsec:Evaluation}, we propose a targeted bait distribution strategy, which takes advantage of heterogeneous mixing, for combating the outbreak.

 \subsection{Homogeneous Mixing Scenario \label{subsec:Homogeneous}}       
 Mosquitoes, like all life forms, have a specific energy requirement, which they can obtain either from taking a blood meal or a sugar meal. However, mosquitoes may prefer blood meals over sugar meals, as the former contain a rich variety of nutrients, such as proteins, which are used by mosquitoes for oviposition. Conversely, if the sugar baits are positioned near mosquito breeding sites, or are augmented by chemical odors which are attractive to the mosquitoes, then they may take more sugar meals than blood meals. We incorporate these scenarios in our model by introducing a parameter $p,\ p \in [0,1]$, which represents the preference of mosquitoes for a blood meal, while $q = 1-p$ represents the preference for a sugar meal. Due to the homogeneity of human hosts, during time interval $[t,t+dt)$, the chance that a mosquito bites a human is given by $apNdt/(pN + qNx)$, and the probability that a mosquito takes a sugar meal in time interval $[t,t+dt)$ is given by $aqNxdt/(pN + qNx)$. This is because sugar baits are assumed to be distributed uniformly throughout the population. If mosquitoes do not have any bias towards sugar baits or humans ($p=.5$), then the chance of taking a sugar meal or a blood meal will simply depend on amount of human hosts and sugar baits. 
 \par
 {A susceptible mosquito ingests the parasite by biting an infected human. The probability of this event is $\frac{acpI_h}{Np + qNx}$.}
  For susceptible mosquitoes we can now write the following transition equation.
  \par
{ \footnotesize
    \begin{align*}
 S_m(t+dt) = & S_m(t) + \underbrace{M\delta dt}_{\textrm{Births}}  - \underbrace{\delta S_m(t)dt}_{\textrm{Deaths}}   -  \underbrace{\frac{a\gamma S_m(t)qNxdt}{Np + qNx}}_{\textrm{Ingest Bacteria}} \\  
 -&  \underbrace{\frac{acpS_m(t)I_h(t) }{Np + qNx}dt}_{Ingest Parasite}   + \underbrace{\theta R_m(t)dt}_{\textrm{Effect wears off}}  + o(dt)    
 \end{align*}
 }
 \par
   {A mosquito in the exposed class can either die, enter the removed class if it ingests the bacteria or become infectious}. It takes $\tau$ duration of time for the mosquito to become infectious. During this time, the mosquito can either die or it can move to the removed compartment by taking a sugar meal. The probability that a mosquito which has gametocytes in its midgut becomes infectious is $(1-\delta dt)^{\tau/dt}(1-a\gamma q xdt/(p+qx))^{\tau /dt}$.  The transition equations for the exposed class are 
   \par
   { \footnotesize
   \begin{align*}
    E_m(t + dt) &= E_m(t) - \underbrace{E_m(t)\delta dt}_{\textrm{Deaths}} + \underbrace{\frac{acp{S}_m(t)I_h(t)}{Np + qNx}dt}_{\textrm{Ingest parasite}}  + o(dt)\\
    &- \underbrace{\frac{a\gamma E_m(t)qNx}{Np + qNx} dt}_{\textrm{Ingest bacteria}} 
     - \underbrace{\frac{acp{S}_m(t-\tau)I_h(t-\tau)e^{- \Lambda \tau}  }{Np + qNx} dt}_{\textrm{Parasites reach salivary gland}}
   \end{align*}
   }
   \par
    where $\Lambda = \delta + \frac{a\gamma qNx }{Np + qNx}$.
Similarly, the transition equation for infectious and  removed classes are
\par
{ \footnotesize
\begin{align*}
I_m(t + dt) &= I_m(t) - \underbrace{I_m(t)\delta dt}_{\textrm{Deaths}} +  \underbrace{\frac{acp{S}_m(t-\tau)I_h(t-\tau)e^{- \Lambda \tau}  }{Np + qNx}dt  }_{\textrm{Parasites reach salivary gland}}  \\ &+  o(dt) \\
R_m(t+dt) &= R_m(t) -\underbrace{R_m(t)\delta dt}_{\textrm{Deaths}} - \underbrace{R_m(t)\theta dt}_{\textrm{Effect wears off}} \\&+ \underbrace{\frac{(S_m(t) + E_m(t))aq\gamma Nx}{Np + qNx}}_{\textrm{Ingest bacteria}} + o(dt)
\end{align*}   
   }
   \par
   Taking limit $dt \to 0$, we get the system of differential equations
   { 
  \begin{align}
 \dot{S}_h &= \mu I_h - \frac{abpS_hI_m}{Np + qNx}  \nonumber  \\
   \dot{I}_h &= -\mu I_h + \frac{abpS_hI_m}{Np + qNx}  \nonumber \\
   \dot{{S}}_m &= - \frac{ a\gamma qNx{S}_m }{Np + qNx} - \frac{acp{S}_mI_h}{Np + qNx} + \delta I_m \nonumber  \\ 
   &+ \delta E_m  + (\delta + \theta)R_m  \nonumber \\
   \dot{E}_m &= -\delta E_m  +  \frac{acp\hat{S}_mI_h}{Np + qNx} - \frac{a\gamma E_mqNx}{Np + qNx}  \nonumber \\
   &-   \frac{acp{S}_m(t-\tau)I_h(t-\tau)e^{- \Lambda \tau} }{Np + qNx} \nonumber \\
 \dot{I}_m & =  \frac{acp{S}_m(t-\tau)I_h(t-\tau)e^{- \Lambda \tau} }{Np + qNx}  -  \delta I_m  \nonumber \\
   \dot{R}_m &= \frac{a\gamma qNx({S}_m + E_m)}{Np + qNx} - (\theta + \delta)R_m  \nonumber
\end{align}
For ease of mathematical analysis we reduce the above system of equations by eliminating $E_m$. This is achieved by adding $E_m$ with $S_m$. Defining $\hat{S}_m = S_m + E_m$, we get 
  \begin{align}
\dot{S}_h &= \mu I_h - \frac{abpS_hI_m}{Np + qNx}  \nonumber  \\
\dot{I}_h &= -\mu I_h + \frac{abpS_hI_m}{Np + qNx}  \nonumber \\
\dot{\hat{S}}_m &= - \frac{acp\hat{S}_m(t-\tau)I_h(t-\tau)e^{- \Lambda \tau} + a\gamma qNx\hat{S}_m }{Np + qNx}  \nonumber \\ & + \delta I_m + (\delta + \theta)R_m  \nonumber \\
\dot{I}_m &= \frac{acp\hat{S}_m(t-\tau)I_h(t-\tau)e^{- \Lambda \tau}}{Np + qNx}  - \delta I_m \nonumber \\
\dot{R}_m &= \frac{a\gamma qNx\hat{S}_m}{Np + qNx} - (\theta + \delta)R_m \label{eqn:homo_delay}
\end{align}
}
In a typical Ross-MacDonald  malaria model with delay \cite{Ruan2008raey}, the system of DDE is cooperative \cite{Ruan2008raey,Smith2008}.  It can be shown, \cite{Smith2008}, that for calculating the stability of equilibrium points of a cooperative system of DDEs  with finite and discrete delays using linear stability analysis, the delays can be ignored and linear stability analysis can be performed on the corresponding ODE. However, the above system of delay differential equations, (\ref{eqn:homo_delay}), is a non-cooperative dynamical system because $\frac{\partial I_m}{\partial R_m} \leq 0$, \cite{Smith2008}. This is specifically due to the introduction of the removed class for mosquitoes. Since system (\ref{eqn:homo_delay}) is not a cooperative system we cannot simply ignore the delays. We use the method of characteristic equation to perform stability analysis of the disease free equilibrium point for calculating the reproductive number.        
\par
  Let $i_h = I_h/N,\ i_m = I_m/M,\ r_m = R_m/M$. Now $r_{df}$ is also the disease free equilibrium proportion of mosquitoes in the removed state. Similarly let $i_{h-end}$ and $i_{m-end}$ be the endemic equilibrium proportion of humans and mosquitoes in the infectious state. 
  Let,
  \par
  { \footnotesize
    \begin{align}
   i_{h-end}&= \frac{(R_0-1)\delta \mu (p+qx)^2  }{(1-r_{df})acpe^{-\Lambda \tau}(abmp + \mu(p+qx))} \nonumber \\
   i_{m-end}&= \frac{(R_0-1)(p+qx)^2\mu \delta}{abmp(1-r_{df})\left(\delta (p+qx)+acpe^{-\Lambda \tau} + \delta a\gamma qx/(\delta + \theta) \right) } \nonumber \\
   r_{df} &= \frac{a\gamma qx}{(p+qx)(\delta + \theta) + a\gamma qx} \nonumber
   \end{align}  
   }
   \par
  where 
  \begin{align}
  R_0 &= \frac{a^2bcmp^2(1-r_{df})e^{-\Lambda \tau}}{\delta \mu(p+qx)^2} \label{eqn:R0}
  \end{align}
In the following theorem  we analyze the stability  of the system of DDEs.
\begin{thm}
 If $R_0 \leq 1$ the system (\ref{eqn:homo_delay}) displays a unique equilibrium point $(i_h = 0, i_m = 0, r_m = r_{df})$ which is stable, and if $R_0>1$ then the system exhibits two equilibrium points: $(i_h = 0, i_m = 0, r_m = r_{df})$ and $(i_h = i_{h-end}, i_m = i_{m-end}, r_m = r_{df}(1-i_{h-end}))$. For $R_0 > 1$ the equilibrium point  $(0,0,r_{df})$ is unstable.
   \end{thm}

 \begin{proof}:
  Since, $S_h + I_h = N$ and $\hat{S}_m + I_m + R_m = M$ we can reduce the the system of five  equations (\ref{eqn:homo_delay}) to three equations. Rewriting the above equations in terms of proportions we get:
     
  \begin{align}
  \dot{i}_h &= -\mu i_h + \frac{abmp(1-i_h)i_m}{p + qx}   \nonumber \\
  \dot{i}_m &= \frac{acp(1-i_m(t-\tau)-r_m(t-\tau))i_h(t-\tau)e^{-\tau \Lambda}}{p + qx}  \nonumber  - \delta i_m  \nonumber  \\
  \dot{r}_m &= \frac{a\gamma qx}{p + qx}(1- i_m - r_m) - (\theta + \delta)r_m  \nonumber
  \end{align}

  Setting $\dot{i}_h, \dot{i}_m $ and $\dot{r}_m$ in the above equations to $0$ we get.
 \begin{align}
& i_h\left (\mu + \frac{abmpi_m}{p+qx} \right ) = \frac{abmp}{p+qx}i_m \nonumber \\ 
& i_m\left(\delta + \frac{acpi_he^{-\Lambda \tau}}{p+qx} \right) =  \frac{acp(1 - r_m)e^{-\Lambda \tau}}{p+qx} i_h \nonumber \\
& r_m = \frac{a\gamma qx}{(p+qx)(\delta + \theta) +a\gamma qx}(1-i_m) \nonumber  
\end{align}
  Put $i_h =0$ and $i_m=0$ for obtaining $r_{df}$. $(0,0,r_{df})$ is an equilibrium point of the system as $i_m=0, i_h=0$ satisfies the above equations. To calculate the other equilibrium point, we solve the three equations simultaneously. After some algebraic simplifications we obtain:
\par
{\small
\begin{align*}
 i_h&= \frac{a^2bcmp^2(1-r_{df})e^{-\Lambda \tau}-\delta \mu (p+qx)^2  }{(1-r_{df})acpe^{-\Lambda \tau}(abmp + \mu(p+qx))} \nonumber \\
 i_m&= \frac{a^2bcmp^2(1-r_{df})e^{-\Lambda \tau}-\delta \mu (p+qx)^2}{abmp(1-r_{df})\left(\delta (p+qx)+acpe^{-\Lambda \tau} + \delta a\gamma qx/(\delta + \theta) \right) } \nonumber \\
r_m&= r_{df}(1-i_m) \nonumber 
\end{align*}
}
\par
Substituting for $R_0$, we get
{\footnotesize
 \begin{align}
 i_h&= \frac{(R_0-1)\delta \mu (p+qx)^2  }{(1-r_{df})acpe^{-\Lambda \tau}(abmp + \mu(p+qx))} \nonumber \\
 i_m&= \frac{(R_0-1)(p+qx)^2 \mu \delta}{abmp(1-r_{df})\left(\delta (p+qx)+acpe^{-\Lambda \tau} + \delta a\gamma qx/(\delta + \theta) \right) } \nonumber \\
r_m&= r_{df}(1-i_m) \nonumber 
 \end{align}
 }
 \par
 If $R_0 \leq 1$ then $i_h=0$ and $i_m=0$. Thus, the system has a unique equilibrium point,  $(0,0,r_{df})$. If $R_0 > 1$ the system has two equilibrium points $(0,0,r_{df})$ and $(i_{h-end},i_{m-end},r_{df}(1-i_{h-end}))$. We now proceed to prove the stability of the equilibrium point using a method similar to one in \cite{Ruan2008raey}.
 
 The system can be linearized at $(0,0,r_{df})$ as $i_m$, $i_h$ and $r_m - r_{df}$ are small resulting in negligible higher order terms. Thus we get, 
 \begin{align}
 \dot{i}_h &= -\mu i_h + \frac{abmpi_m}{p + qx}   \nonumber \\ 
 \dot{i}_m &= \frac{acp(1-r_{df})i_h(t-\tau)e^{- \Lambda  \tau}}{p + qx}  - \delta i_m \label{eqn:linearized} \\ 
 \dot{r}_m &= \frac{a\gamma qx}{p + qx}(1-i_m-r_m) - (\theta + \delta)r_m \nonumber  
 \end{align} 
 The characteristic equation of the above system  is given by:
  \begin{align}
 &\left(\lambda + \delta + \theta + \frac{a\gamma qx}{p+qx}\right) \nonumber \\
 & \times \left(\lambda^2 + \lambda(\delta + \mu) + \delta \mu - R_0e^{-\lambda \tau}\delta \mu \right) = 0 \label{eqn:characteristic}
 \end{align}
 
 The system is stable at the given point if and only if all the real parts of roots of the above equation are negative. Let   \begin{align}F(\lambda,\tau) = \lambda^2 + \lambda(\delta + \mu) + \delta \mu - R_0e^{-\lambda \tau}\delta \mu = 0 \label{eqn:transcend}\end{align}  The system is stable if and only if real parts of roots of $F(\lambda,\tau)$ are negative. We now consider three different cases for exploring the roots of $F(\lambda, \tau)$.
  (i)$R_0<1$: We first show the existence of roots for $F(\lambda,\tau) = 0$. Now, $F(\lambda,\tau) = 0$ has two roots and since $R_0<1$ both are negative. The two roots are given by
   \begin{align*}
  \lambda_\pm = \frac{-(\mu + \delta) \pm \sqrt{(\delta+\mu)^2 - 4\delta \mu (1-R_0)}}{2}
  \end{align*}
 The two roots cannot be the same, as that would require $(\delta-\mu)^2 =  - 4R_0$.  Also,
   \begin{align*}
   F_{\lambda}^{'}(\lambda, 0) = 2\lambda +  (\delta+\mu)
  \end{align*}
 Thus, $F_{\lambda}^{'}(\lambda_\pm, 0) \neq 0$. From the implicit function theorem and continuity of $F(\lambda,\tau)$, the equation (\ref{eqn:transcend}) has roots for all $\tau > 0$. We now show that the real part of these roots must be negative.
 \par
 $F(0,\tau)>0$ and $F_\lambda^{'}(\lambda,\tau) > 0$ for all $\lambda \geq 0, \ \tau >  0$. Thus equation (\ref{eqn:transcend}) has no zero roots and positive roots for all positive $\tau$. We claim that equation (\ref{eqn:transcend}) has no purely imaginary roots. Suppose it has a pair of imaginary roots $\pm \omega i$ for some $\tau$. Multiplying the equation obtained by substituting $\omega i$ in equation (\ref{eqn:transcend}) with its conjugate we get.

 \begin{align*}
 \omega ^4 + (\mu^2 + \delta^2)\omega^2 + \delta^2\mu^2 - (R_0 \delta \mu)^2 = 0  
 \end{align*}  

 Now, $\omega$ must be a positive root of the above equation, but since $R_0< 1$ the above equation cannot have non-negative real roots, which is a contradiction. Thus, (\ref{eqn:transcend}) has no purely imaginary roots and no positive roots, and hence real part of the roots must be negative. 
 
 (ii)$R_0 = 1$: then $F(0,\tau) = 0$ and $F_\lambda^{'}(\lambda,\tau) > 0$ for $\lambda \geq 0, \ \tau >  0$. Equation (\ref{eqn:transcend}) has zero root and no positive root. Using an argument similar to case (i) it can be shown that the root other than $\lambda = 0$ has a negative real part. 
 
 (iii)$R_0>1$: then $F(0,\tau) < 0$ and $F_\lambda^{'}(\lambda,\tau) > 0$ for all  $\lambda \geq 0, \ \tau > 0$.  Therefore equation (\ref{eqn:transcend}) must have a unique positive real root for all $\tau > 0$. Hence, the point $(0,0,r_m^*)$ is unstable if $R_0 > 1$.
 \par
 This completes the proof.
 \end{proof}   
\par
 An outbreak is possible only when $R_0 > 1$. The factor $(qx)^2$ can be seen in the denominator of the $R_0$ expression. Thus, for low $p$, $R_0$ is inversely proportional to bait density ($x$) times attractiveness of the baits $q$. We term this product $xq$ as effective bait density.
 \par
  Notice that substituting $x=0$ (no paratransgenesis) results in $R_0 = a^2bcme^{-\delta \tau} / \mu \delta$, which is the reproductive number obtained from the original Ross-MacDonald  model \cite{Smith2004a,Smith2012}.  Simulations show that for $R_0>1$ case, the system eventually reaches a  stable equilibrium point. This suggests that the equilibrium point $(i_{h-end},i_{m-end},r_{df}(1-i_{h-end}))$ must be stable when $R_0>1$.  
 
\subsection{Heterogeneous Mixing Scenario \label{subsec:Heterogeneous}}
In the previous section we assumed that mosquitoes do not differentiate between humans for obtaining a blood meal. However, many studies \cite{Smith2004b,Smith2005,Bejon2010,Bousema2012} suggest that this assumption may not always hold. Empirical studies \cite{Woolhouse1997,Smith2007} have suggested that malaria infections follow a statistical pattern known as the $20/80$ rule: $20 \%$ of the humans are responsible for $80 \%$ of the infections. This indicates that mosquitoes favor a certain class of individuals over others for obtaining a blood meal. According to a recent study \cite{Verhulst2011}, this `attractiveness' may also be dependent on the nature of the skin microbiota.  Thus, the heterogeneity in number of mosquito bites may arise from not only the environmental factors such as proximity of a settlement to a swamp but also from biological factors such as skin microbiota. For example populations located near a swamp (hot spot) will suffer from a higher bite rate than those located farther away from the swamp. For settlements which are equidistant to the swamp, biological factors like skin microbiota may play a role. The resulting heterogeneity can be captured through the probability distribution of the attractiveness over the human population. 
\par
We incorporate this heterogeneity in the model by classifying the human population on the basis of their `attractiveness' to mosquitoes. Humans belonging to class $i$ are assumed to have an attractiveness factor $k_i$. Let $T$ be the total number of classes. Let $P(i)$ be the probability that a randomly chosen human belongs to class $i$, and $N_i$ be the number of human hosts in class $i$, thus $P(i) \approx N_i/N$. Therefore the chance that a mosquito bites an individual belonging to class $i$ in time interval $[t,t+dt)$ is given by $ak_ipN_idt/(p \sum_{i}^{T}k_i N_i + qNx) $; while the chance that it takes a sugar meal in time interval $[t,t+dt)$ is given by $aqNxdt/(p \sum_{i}k_i N_i + qNx)$. Let $I_h^i$ be the number of infected human hosts in class $i$ and let $i_h^i = I_h^i/N_i$. Let $\bar{k}$ denote the mean attractiveness factor among the human population $\bar{k} = \sum_{i}^{T}k_iP(i)$, and let $\hat{k}$ denote the second moment, i.e., $\hat{k} = \sum_{i}^{T}k_i^2P(i)$. Application of mean field approximation results and combining the exposed and susceptible compartments, we obtain the following set of DDEs:
\par
{\small 
 \begin{align} 
 \dot{S}_h^i &= -\frac{abk_ip(N_i - I_h^i)I_m}{p \sum_{i}k_i N_i + qNx} + \mu I_h^i  \nonumber \\
\dot{I}_h^i &= \frac{abk_ip(N_i - I_h^i)I_m}{p \sum_{i}k_i N_i + qNx} - \mu I_h^i \label{eqn:I_h^i hetero}  \\
\dot{I}_m &= \frac{acp(M-I_m(t-\tau) -R_m(t-\tau))e^{-\zeta \tau}\sum_i k_i I_h^i(t-\tau)}{p \sum_{i}k_i N_i + qNx} \nonumber \\& -\delta I_m  \nonumber \\ 
\dot{\hat{S}}_m &= -\frac{acp(M-I_m(t-\tau) -R_m(t-\tau))e^{-\zeta \tau}\sum_i k_i I_h^i(t-\tau)}{p \sum_{i}k_i N_i + qNx} \nonumber \\ 
&-\frac{a \gamma q Nx(M-I_m-R_m)}{p \sum_{i}k_iN_i + qNx} +\delta I_m + (\delta + \theta)R_m \nonumber \\ 
\dot{R}_m &= \frac{a \gamma q Nx(M-I_m-R_m)}{p \sum_{i}k_iN_i + qNx} - (\theta + \delta)R_m  \nonumber
\end{align}   
}
\par
 where $\zeta = \delta + a\gamma q x/(p\bar{k}+qx)  $. Observe that $\sum_{i}k_i N_i / N = \sum_{i}k_i P(i) = \bar{k} $. Let  
 \begin{align*}
 \phi := \sum_i k_i I_h^i/N
 \end{align*}
  We obtain $\dot{\phi}$ from (\ref{eqn:I_h^i hetero}).{ Also, $\dot{i}_h = \frac{1}{N}\sum_{i=1}^{N}i_h^i$.}  The above system can be written as:
  \begin{align} 
  \dot{i}_h &= \frac{abmp(\bar{k} - \phi)i_m}{p\bar{k} + qx} - \mu i_h \nonumber \\ 
  \dot{i}_m &= \frac{acp(1-i_m(t-\tau) - r_m(t-\tau))e^{-\zeta \tau}\phi(t-\tau)}{p\bar{k} + qx} \nonumber \\ &-\delta i_m \nonumber \\ 
  \dot{r}_m &= \frac{a\gamma q x}{p\bar{k} + qx}(1-i_m-r_m) - (\delta + \theta)r_m \nonumber \\
  \dot{\phi}&= \frac{i_m(\hat{k}-\sum_i k_i^2 i_h^i)abmp}{p\bar{k}+ qx}  - \mu \phi \label{eqn:DDE_hetero_2}
  \end{align}     
  This is potentially a large system of DDEs, since the number of classes, $T$, can be very high. We now show than even for such a large system of DDEs, the stability analysis of the disease free equilibrium state can be carried out and a closed form expression for $R_0$ can be obtained.
 \par  
 Let $r_{df}$ be the disease free equilibrium proportion of mosquitoes in the removed state. 
 \begin{align*}
 r_{df} = \frac{a\gamma qx}{ (p\bar{k}+qx)(\delta + \theta) + a\gamma qx} 
 \end{align*}
 We define $R_0$ for the heterogeneous case, 
 \begin{align}
 R_0 &= \frac{a^2bcmp^2 \hat{k}(1-r_{df})e^{-\zeta \tau}}{\delta \mu (\bar{k}p+qx)^2}  \label{eqn:R_0_hetero} \\
  \end{align}
 We now discuss the stability of the above system in the following theorem. 
\begin{thm}\label{thm:hetero}
The equilibrium point $(0,0,r_{df},0)$ of system (\ref{eqn:DDE_hetero_2}) is stable if and only if $R_0 \leq 1$.  If any equilibrium point ($i_h,i_m,r_m,\phi$) other than $(0,0,r_{df},0)$ exists then such a point must satisfy $i_h > 0$, $i_m >  0$. 
\end{thm}
 \begin{proof}:
  It is easy to see that $(0,0,r_{df},0)$ is an equilibrium point of system (\ref{eqn:R_0_hetero}).   Thus the equilibrium point $(0,0,r_{df},0)$ of system (\ref{eqn:DDE_hetero_2}) is stable if and only if $R_0 \leq 1$.
 We can linearize the system at $(0,0,r_{df},0)$ since $i_m$, $i_h$, $\phi$ and $r_m - r_{df}$ are very small and can be neglected. The linearized system is given by    
  \begin{align} 
  \dot{i}_h &= \frac{abmp\bar{k}i_m}{p\bar{k} + qx} - \mu i_h \nonumber \\ 
  \dot{i}_m &= \frac{acp(1-r_{df})e^{-\zeta \tau}\phi(t-\tau)}{p\bar{k} + qx}  -\delta i_m \nonumber \\ 
  \dot{r}_m &= \frac{a\gamma q x}{p\bar{k} + qx}(1-i_m-r_m) - (\delta + \theta)r_m \nonumber \\
  \dot{\phi}&= \frac{abmp\hat{k}i_m}{p\bar{k}+ qx} - \mu \phi \label{eqn:DDE_hetero_linearzed}
  \end{align}     
The characteristic equation of the DDEs is given by  
  \begin{align*}
&  (\lambda+\mu)\left(\lambda + \delta + \theta + \frac{a\gamma qx}{p\bar{k}+qx}\right) \\ & \times \left(\lambda^2 + \lambda(\delta + \mu) + \delta \mu - R_0e^{-\lambda \tau}\delta \mu \right) = 0
  \end{align*}
  \par
 This is very similar to equation (\ref{eqn:characteristic}), thus the proof for stability can be completed using the same arguments. 
\par
  Now, if there is any equilibrium point other than $(0,0,r_{df},0)$, then it must satisfy $i_h^* > 0$ and $i_m^* >  0$. Let us assume that the converse is true. Then we have three cases: (i) $i_h^* = 0 , i_m^* = 0$; (ii) $i_h^* > 0, i_m^* =  0$; (iii) $i_m^* > 0, i_h^* = 0$.  
  \par
  If case (i) is true then $\phi^* = 0$ and $r_m^* = r_{df}$ which is the same as equilibrium point $(0,0,r_{df},0)$. If case (ii) is true, we substitute $i_m = 0 $ in system (\ref{eqn:DDE_hetero_2}). We get,  $\dot{i}_h = 0$ if and only if $i_h^* = 0 $ or $\mu = 0$. Since $\mu > 0$, $i_h^*$ must be $0$. If case (iii) is true then $\phi^* = 0$, we substitute $\phi = 0$ in system (\ref{eqn:DDE_hetero_2}). We get, $\dot{i}_m = 0$ if and only if $i_m^* = 0$ or $\delta = 0$. Since $\delta > 0$, $i_m^*$ must be $0$. This completes the proof.
  \end{proof}
  \par
   The threshold quantity $R_0$ is an increasing function of $\hat{k}$. If $X,Y$ are two non negative discrete random variables with the same mean, and if $X \succeq Y$, then $E[X^2] \geq E[Y^2]$. Thus, heavy (fat) tailed distributions will have a higher $R_0$ than light tailed distributions. Also, similar to the homogeneous case, $R_0$ is inversely proportional to the effective bait density for small values $p\bar{k}$.
  \par
 { The homogeneous mixing case is a specific instance of the heterogeneous case. We can recover the $R_0$ for the homogeneous case from the heterogeneous case as follows. In the homogeneous case all humans have the same attractiveness factor; let that be $k$. Hence, $\bar{k} = k$ and $\hat{k}= k^2$. If $k = 1$, then we obtain $R_0 = \frac{a^2bcmp^2 (1-r_{df})e^{-\Lambda \tau}}{\delta \mu (p+qx)^2}$ which is same as the $R_0$ for the homogeneous case. }   

\subsection{Evaluation of a Targeted Bait Distribution Strategy \label{subsec:Evaluation}}   
Malaria `hotspots' \cite{Bousema2012}, i.e., geographical areas where malaria transmission is much higher than the average, give rise to heterogeneous biting. The presence of such hotspots may result in a higher reproductive number which can sustain the outbreak. A study carried out in Kenya, \cite{Ernst2006}, discovered that the chance of catching malaria inside the hotspot was $2.6 - 3.2$ times more than that outside the hotspot. However, in order to target these hotspots, one has to first detect these hotspots. Hotspots can be detected by analyzing data obtained from : asymptomatic parasite carriage \cite{Bejon2010}, serological testing on antigens \cite{Bousema2010b, Jacklin2014}, mosquito density \cite{Bousema2010}, exposure to infected mosquitoes \cite{Bousema2010}. In this subsection we propose  a mechanism of selectively targeting population with high attractiveness factor.         
\par
 Instead of uniformly distributing the sugar-baits, we propose a targeted distribution of sugar baits, i.e., let $B_i$ be the number of sugar baits distributed around or inside the residences of individuals belonging to class $i$. Now, since these sugar baits are distributed in the close vicinity of individuals with attractiveness factor $k_i$, the attractiveness factor for these  baits is assumed to be approximately $k_i$. Hence, mosquitoes have a preference $qk_i$ for baits belonging to class $i$. The prime question of interest is, that given a limited amount of sugar baits $B = Nx$ what is the amount of sugar baits must be deployed around individuals of class $i$, i.e., $B_i$ for minimizing the incidence of malaria.
\par
Since these baits are deployed in the vicinity of humans, there will be constraints on the number of deployed baits due to spatial and ethical reasons. Let $C_{i}$ be the maximum number of baits that can be deployed around individuals in class $i$. If $\sum_{i}C_{i} \leq Nx$, then the question of allocating baits does not arise, as all (or a fraction) of the available baits will be deployed. If $\sum_{i}C_{i} > Nx$ then we need to find an efficient distribution scheme which minimizes the reproductive number. The reproductive number for targeted baits can be obtained from (\ref{eqn:R_0_hetero}). 
 
\begin{align}
R_0 &= \frac{a^2bcmp^2 \hat{k}(1-r_{df})e^{-\zeta \tau }}{\delta \mu (p\bar{k}+qy)^2} \label{eqn:R0_hetero_2} \\
y   &= \sum_i k_iB_{i}/N\nonumber \\
\zeta  &= \delta + a\gamma qy/(p\bar{k}+qy) \nonumber 
\end{align} 
 
We formulate an optimization problem for calculating a distribution scheme which minimizes the reproductive number.
 
   \begin{align*}
   \begin{aligned}
   & \underset{B_{i}}{\text{minimize}}
     \ \ R_0 \\
   & \text{subject to} \ \ \\
   &		B_{i} \leq C_{i} \ \forall \ i = 1, \ 2, \dots T \\
   &		\sum_i B_{i} = Nx \\
   &		B_{i} \geq 0 \ \forall \ i = 1, \ 2, \dots T
    \end{aligned}
   \end{align*}    
 Allocating highest number of baits to the most attractive class subject to the constraints seems a natural strategy. Below we show mathematically that such a strategy is also globally optimal.
\begin{thm}
The following strategy achieves maximum: Sort the classes starting from lowest attractiveness factor to the highest, i.e., $k^{(1)} < k^{(2)} < \dots < k^{(T)}$ and let $B_{i}$ and $C_{i}$ be the number of baits and the constraint for the ordered class $i$. First allocate $C_{T}$ baits to class $T$ then $C_{T-1}$ baits to class $T-1$ and so on till the available baits are exhausted. The solution is found to be a unique global maximizer.
\end{thm}
 \begin{proof}:
 It can be shown that the objective function $R_0$ is strictly decreasing with $y$. Thus, maximizing $y$ or $Ny$ will minimize $R_0$. Dividing $B_{i}$ by $B = Nx$ and rewriting:
 
  \begin{align*}
    \begin{aligned}
    & \underset{B_{i}}{\text{maximize}}
      \ \ \sum_i k_iB_{i}/B \\
    & \text{subject to} \ \ \\
    &		B_{i} \leq C_{i} \ \forall \ i = 1, \ 2, \dots T \\
    &		\sum_i B_{i}/B = 1 \\
    &		B_{i} \geq 0 \ \forall \ i = 1, \ 2, \dots T
     \end{aligned}
    \end{align*}  
We use stochastic dominance for the proof. { A random variable $X$ stochastically dominates another random variable $Y$ ($X \succeq Y$) if and only if $\mathbb{P}(X>x) \geq \mathbb{P}(Y>x), \ \forall \ x \ \in \ \mathcal{R}$}. Notice, that the objective function is a convex combination. We define $X$ to be a random variable with probability mass function given by $\mathbb{P}(X =i) = B_{i}^*/B$ where $B_{i}^*$ is the proposed allocation. We define $Y$ as a random variable with probability mass function given by $ \mathbb{P}(Y = i ) = B_{i}^a/B $ where $ B_{i}^a $ is any other allocation.
 By construction, $\mathbb{P}(X>x) \geq \mathbb{P}(Y>x) ,\  \forall \  x = 1,\ 2 \dots T$. Hence, $X \succeq Y $  and therefore $E[f(X)] \geq E[f(Y)]$ for a non decreasing function $f$. If $f(i) = k^{(i)} $ then $E[f(X)]$ and $E[f(Y)]$ are the objective functions and hence the proposed allocation is optimal. Furthermore, it is unique as $k^{(i)} \neq k^{(j)}$ for $i \neq j$.  This completes the proof.
 \end{proof}
 
 \par
 { The optimization problem assumes that bait would have attractiveness factor similar to the individuals around which it is deployed. However, baits are deployed near houses and individuals with similar attractiveness factor may not live in the same house. Thus, the attractiveness factor of the bait may no longer be the same as that of one particular individual living in the house. One way of solving this is to use the average attractiveness factor of the house (averaging across individuals in the same house) as the attractiveness factor of the bait. The solution would be similar to the above problem: allocate the highest number of possible baits to the house with the largest attractiveness factor, then do the same for the house with the second highest attractiveness factor and so on.  Since the above optimization problem uses a more fine grained approach (since it is based on attractiveness of individuals rather than a group of individuals with dissimilar attractiveness factor), this averaging would result in a sub-optimal reduction of $R_0$. }

\section{Numerical results \label{sec:Numerical}}

   \begin{table}[t]
        \centering
  \caption{Parameter Values.}
  \label{table:Param Table}
  \begin{tabular}{l l l l}
\hline\noalign{\smallskip}
  Symbol   &  Values  &  Unit  &  References \\
\noalign{\smallskip}\hline\noalign{\smallskip}
  $a$      & 0.5 & bites per day & \cite{Menach2005}  \\
  $b$      & 0.5 & n/a  &  \cite{Smith2004a,Menach2005} \\ 
  $c$      & 0.5  & n/a &  \cite{Smith2004a,Ruan2008raey} \\
  $m$      & 2 & n/a & \cite{Ruan2008raey} \\
  $\tau$ & 10 & days &   \cite{Menach2005,Hancock2011a}  \\
  $\gamma$ & 0.8 & n/a & \cite{Wang2012} \\
  $1/\theta$ & 5 & days  & \cite{Wang2012} \\
  $\mu$    & 0.05 & recoveries per day  & \cite{Ruan2008raey}  \\
  $\delta$ & 0.1  & births/deaths per day   & \cite{Smith2004a,Menach2005} \\
\noalign{\smallskip}\hline
  \end{tabular}
        \end{table}
    \begin{figure}
    \centering
      \includegraphics[width=0.44 \textwidth]{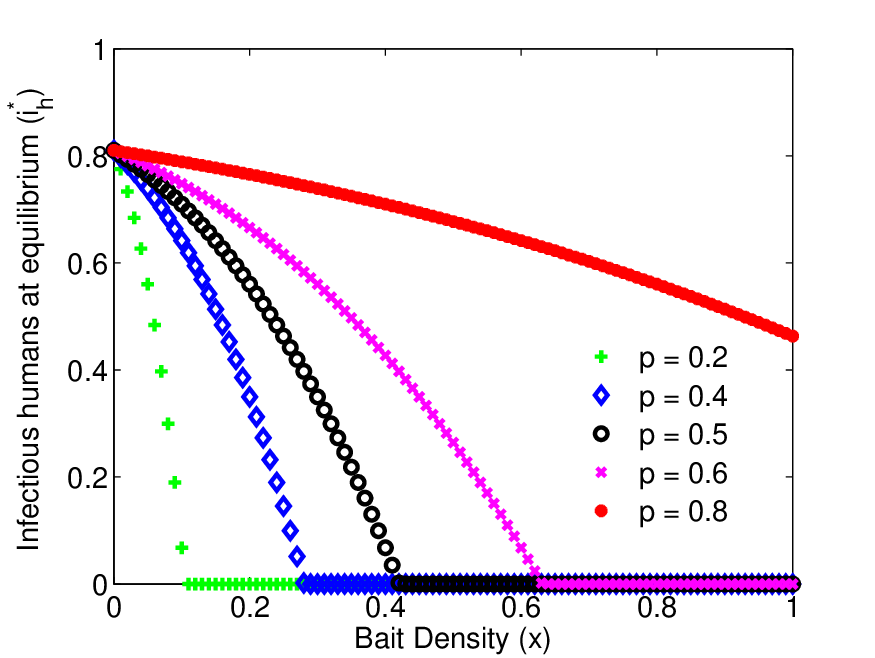}
    \caption{Proportion of infected humans at equilibrium versus bait density for various $p$. See Table \ref{table:Param Table} for parameter values}
    \label{fig:Ih-p}
 \end{figure}
       \begin{figure}
       \centering
    \includegraphics[width=0.44\textwidth]{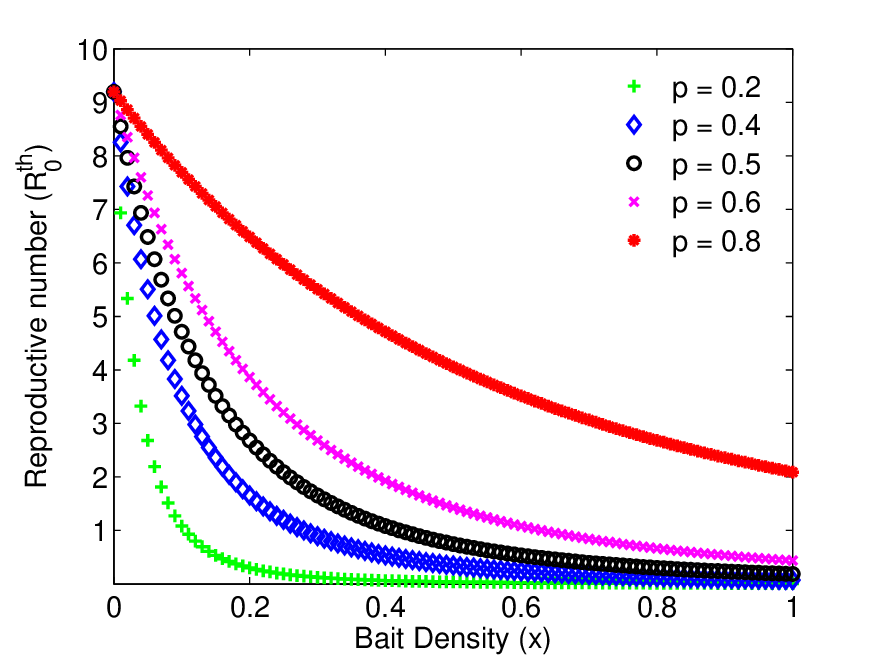}
    \caption{Reproductive number versus bait density for various $p$. See Table \ref{table:Param Table} for parameter values.}
    \label{fig:Reproductive_Num_Homo}
    \end{figure}
In order to answer the question:  would  paratransgenesis be useful if implemented, realistic parameter values are required, and a literature survey was performed for { obtaining parameter values}. 
\par
It takes about $20-100$ days for humans to recover from malaria \cite{Smith2004a,Menach2005,Ruan2008raey},  we assume that the recovery rate for humans $\mu = 0.05$.  The biting rate of mosquitoes on humans is about $0.3-0.9$ per day \cite{Smith2004a,Menach2005,Kileen2007}; here we take $a = 0.5$ per day \cite{Menach2005}. The proportion of bites by infectious mosquitoes which lead to disease, $b$, and the proportion of bites by susceptible mosquitoes on infected humans which causes mosquitoes to become exposed, $c$, are taken to be $0.5$, \cite{Smith2004a,Menach2005,Ruan2008raey}. The mosquito density $m$, or the number of mosquitoes per human, is taken to be $2$ \cite{Ruan2008raey}. The rate of emergence of adult mosquitoes, and their death rate, $\delta$, is about $0.1-0.2$ per day \cite{Smith2004a,Menach2005,Hancock2011a};  we take $\delta = 0.1$ per day. The incubation period of the parasite within the mosquitoes depends on environmental factors like temperature and humidity and it varies with species, \cite{Beier1998}. We assume an incubation period, $\tau$, of $10$ days, \cite{Menach2005,Hancock2011a}. The genetically modified bacteria are effective $80-84\%$ of the time, \cite{Wang2012}; we assume that efficacy of paratransgenesis, $\gamma = 0.8$.  { Since the measurement was done only in the lab, $\gamma = 0.8$ may not hold in the field. Therefore, we also perform numerical analysis for a varying $\gamma$}.  The parasite killing effect of genetically modified bacteria inside the mosquito's midgut lasts for at least $4$ days \cite{Wang2012}. We assume $\theta$ to be $0.2$ per day.  These parameter values are summarized in Table \ref{table:Param Table}.
\par
Since the data concerning the distribution of the attractiveness factor was not available with us, we only perform numerical analysis for the homogeneous mixing scenario. However, given the data, one needs to simply calculate the first, $\bar{k}$, and second moment, $\hat{k}$ of the distribution for computing the reproductive number. Note that the following figures are obtained by numerically evaluating the analytical expressions for $R_0,$ and endemic equilibrium point $i_h$ ; they are not a product of simulations.
        \begin{figure}
        \centering
     \includegraphics[width=0.44\textwidth]{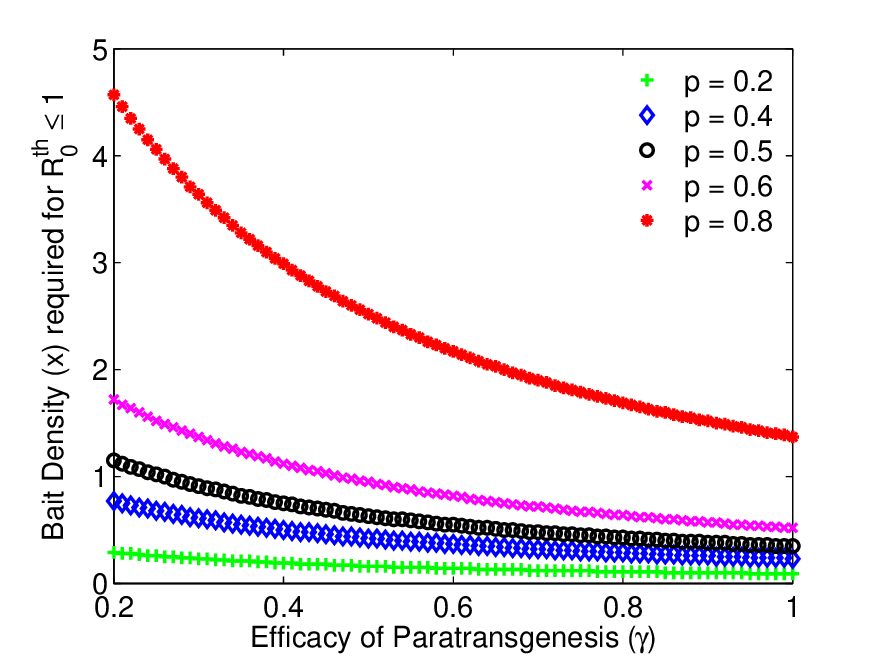}
     \caption{Bait density required for $R_0 \leq 1$ versus $\gamma$ for various $p$. See Table \ref{table:Param Table} for parameter values.}
     \label{fig:Reproductive_Num_gamma}
     \end{figure}

Fig.  \ref{fig:Reproductive_Num_Homo} shows the proportion of baits required for a given $R_0$ for different values of $p$ and parameters listed in Table \ref{table:Param Table}. For $p = 0.2$, approximately one bait per $10$ individuals would be required for keeping $R_0 < 1$, while for $p = 0.8$ the number would be significantly more than the population size. Fig.  \ref{fig:Ih-p} shows the corresponding proportion of infected humans at equilibrium.

 { To further study the usefulness of paratransgenesis we plot the bait density required to drive $R_0$ to $1$ for varying $\gamma$ and $p$ in Fig.  \ref{fig:Reproductive_Num_gamma}. Clearly, a low $p$ and high $\gamma$ is highly desirable. More importantly, the figure suggests that a low value of $p$ is more important than a high value of $\gamma$.  Even if $\gamma$ is small, a low value of $p$ can ensure that the bait density does not cross $1$.  }
        \begin{figure}
        \centering
     \includegraphics[width=0.44\textwidth]{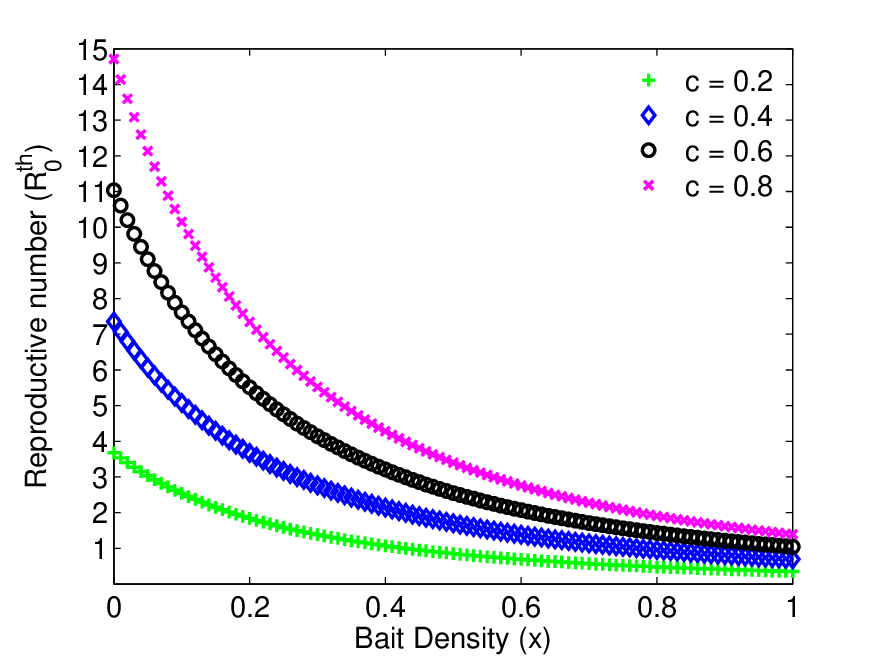}
     \caption{Reproductive number versus bait density for various $c$, $\gamma = 0.3, \ p = 0.5$}
     \label{fig:Reproductive_Num_C}
     \end{figure}
            \begin{figure}
            \centering
         \includegraphics[width=0.44\textwidth]{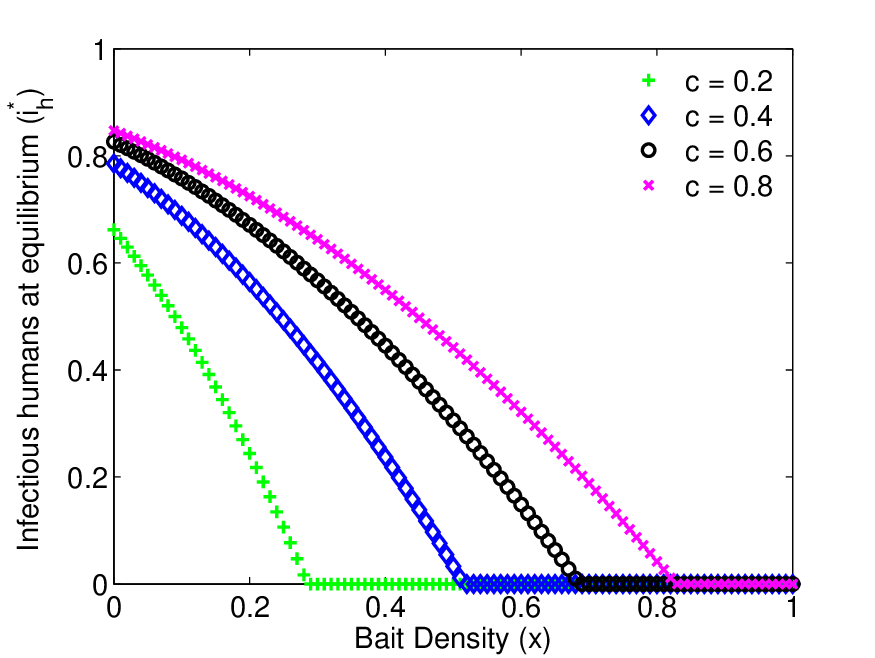}
         \caption{Proportion of infected humans at equilibrium versus bait density for various $c$, $\gamma = 0.3, \ p = 0.5$.}
         \label{fig:Ih_c}
         \end{figure}
 
 { In Fig. \ref{fig:Reproductive_Num_C} we plot $R_0$ versus bait density for varying values of an environment parameter $c$. Fig. \ref{fig:Ih_c} provides the corresponding values of $i_h^*$. Recollect that $c$ is the probability that mosquito ingests the parasite after biting an infected individual. Thus, a low value of $c$ is desirable. As shown by Fig. \ref{fig:Ih_c}, even for a high value of $c$ and low $\gamma$, paratransgenesis can be very helpful in reducing the number of infected humans.}
       
\section{Discussion \label{sec:Discussion}}
Our analyses demonstrate that paratransgenesis can be a viable strategy in controlling malaria under certain conditions. In the homogeneous mixing case, although Fig. \ref{fig:Reproductive_Num_Homo} shows that increase in bait density is associated with rapid decline of the reproductive number, for $p=0.8$ the amount of baits required to keep $R_0$ under unity becomes a multiple of the population size which may make adoption of this strategy impractical. In scenarios where mosquitoes prefer sugar meal over a blood meal ($p < 0.5$) the bait density required to keep $R_0<1$ is less than the population size, paratransgenesis can be viable. This is also evident from the analytical expression of $R_0$ for both homogeneous and heterogeneous case. For a low $p$, $R_0$ is inversely proportional to effective bait density.
\par
 { From Fig. \ref{fig:Reproductive_Num_gamma} one can conclude that a low value of $p$ is more important than a high value of $\gamma$. Thus, efforts directed at making sugar baits more attractive to mosquitoes (reducing $p$) may yield better results than efforts at making paratransgensis more efficient (increasing $\gamma$).  }
\par
The usefulness of the targeted bait allocation strategy depends on the constraints. In scenarios where the constraints are too restrictive, targeted bait distribution may have little effect. { For example, if the constraints dictate that baits cannot be added near populations with high attractiveness factor (hot spots) then the targeted strategy may not be useful.} Another interesting observation is that the effectiveness of paratransgenesis in lowering the reproductive number increases with the incubation time $\tau$. This is because the expression of the reproductive number, equations (\ref{eqn:R0}) and (\ref{eqn:R_0_hetero}) contain the term $e^{-\Lambda \tau}$ and $e^{-\zeta \tau}$ respectively. 
\par
 We did not consider the effect of incubation time of plasmodium inside humans, but our analysis can be easily extended to incorporate that factor. Analyses in \cite{Ruan2008raey} suggest that the effect of incorporating this factor would result in the reduction of the reproductive number, $R_0$,  by a factor of  $e^{-\delta \tau_h}$ where $\tau_h > 0$ is the incubation time within humans. 
\par
It should be noted that malaria control using paratransgenesis is not a vector control strategy, but a disease control strategy, and hence it may not reduce the inconvenience caused by mosquito bites. However, it can be used along with traditional vector control measure such as ITNs \cite{Wang2012} to protect against bites. Thus, paratransgenesis combined with odor augmented baits, ITNs  and targeted bait allocation can prove to be a very effective malaria control strategy.  Our results, in particular the proposed targeted bait distribution strategy, have implications on malaria control policies. The targeted approach takes into account heterogeneous biting and provides an extremely efficient way of allocating resources. Furthermore, it can be implemented in a community by first gathering data on the biting pattern of mosquitoes \cite{Bousema2012} and then placing appropriate number of baits around or inside residential places. Malaria is a global issue, and we believe that our results contribute towards its control and eradication.

\bibliographystyle{elsarticle-num}
\bibliography{Malaria_Paratransgenesis}   

\end{document}